\input{style/arxiv-general.cfg}
\documentclass[MSNbibl,number,citesort,seceqn,dvips]{arxbj}
\makeatletter
   \@ifpackageloaded{graphicx}{}{\usepackage{graphicx}}
\makeatother
\usepackage{upgreek,mathrsfs}
\aid{0}
\volume{21}
\issue{2}
\pubyear{2015}
\firstpage{760}
\lastpage{780}
\doi{10.3150/13-BEJ586} % kopijuoti is 'New paper accepted'

\makeatletter

\newcommand{\rright}{\right}
\newcommand{\lleft}{\left}
\newtheorem{prop}{Proposition}[section]
\newtheorem{cor}[prop]{Corollary}
\newtheorem{lemme}[prop]{Lemma}
\newremark{Remark}{Remark}
\newtheorem{them}[prop]{Theorem}
\newtheorem{theorem}[prop]{Theorem}

\newcommand{\N}{\mathbb{N}}
\newcommand{\Z}{\mathbb{Z}}

\newcommand{\R}{\mathbb{R}}

\def\1{{\mathbf{1}}}
\renewcommand{\pi}{\uppi}

\def\sfrac#1#2{#1/#2}
\def\vfrac#1#2{(#1)/#2}

\def\sklfrac#1#2{(#1/#2)}
\makeatother

\begin{document}
\begin{frontmatter}

\title{Fluctuations of the power variation of fractional Brownian
motion in Brownian time}
\runtitle{Fluctuations of the power variation of fractional Brownian
motion in Brownian time\hspace*{10pt}}

\begin{aug}
%%%% inicialai - be tarpu
\author[A]{\inits{R.}\fnms{Raghid} \snm{Zeineddine}\corref{}\ead[label=e1]{raghid.zeineddine@univ-lorraine.fr}}% \and
%%\runauthor{} %% auto
\address[A]{Institut Elie Cartan de Lorraine, Universit\'e de Lorraine,
campus Aiguillettes, 54506 Vandoeuvre L\`{e}s Nancy Cedex, France.
\printead{e1}}
\end{aug}

% HISTORY:
\received{\smonth{4} \syear{2013}}
\revised{\smonth{12} \syear{2013}}

% ABSTRACT
%
\begin{abstract}
We study the fluctuations of the power variation of fractional Brownian
motion in Brownian time.
\end{abstract}

% KEYWORDS
% visi is mazosios raides ir pagal abecele
%

\end{frontmatter}

%s1 #&#
\section{Introduction}\label{sec1}

Studying the variations of a stochastic process is of fundamental
importance in probability theory.
In this paper, we are interested in the fractional Brownian motion in
Brownian time, which is defined as follows.
Consider a fractional Brownian motion $X$ on $\R$ with Hurst parameter
$H\in(0,1)$, as well as
a standard Brownian motion $Y$ on $\R_+$ independent from $X$. The
process $Z=X\circ Y$ is the so-called
\textit{fractional Brownian motion in Brownian time} (F.B.M.B.T. in short).
It is a self-similar process (of order $H/2$) with stationary
increments, which is not Gaussian.
When $H=1/2$, one recovers the celebrated iterated Brownian motion.

In recent years, starting with the articles of Burdzy \cite{1,2}, there
has been an increased interest in
iterated processes in which one changes the time parameter with
one-dimensional Brownian motion, see, for example, \cite{3,4,5,6,7} to
cite but a few. In the present paper, we are concerned with the study
of the fluctuations of the $p$th variation of $Z$ for any integer $p$,
defined as
\[
R_n^{(p)}(t)=\sum_{k=0}^{\lfloor2^n t\rfloor-1}
(Z_{(k+1)2^{-n}}- Z_{k2^{-n}})^p,\qquad  n\in\mathbb{N}, t\geq0.
\]
At this stage, it is worthwhile noting that we are dealing with the
$p$th variations of $Z$ in the classical sense when $p$ is even
whereas, when $p$ is odd, we are rather dealing with the \textit
{signed} $p$th variations of $Z$.
The interested reader may read \cite{13,14} in order to find relevant
information about power variations.

After proper normalization, we may expect the f.d.d. convergence to a
non-degenerate limit (to be determined) of
\[
S_n^{(p)}(t)=2^{-n\kappa}\sum
_{k=0}^{\lfloor2^n t\rfloor-1} \bigl((Z_{(k+1)2^{-n}}-
Z_{k2^{-n}})^p-E\bigl[(Z_{(k+1)2^{-n}}- Z_{k2^{-n}})^p
\bigr] \bigr),\qquad  n\in\mathbb{N}, t\geq0,
\]
for some $\kappa>0$ to be discovered.
To reach this goal, a classical strategy consists in expanding the
power function $x^p$ in terms of Hermite polynomials. Doing so, our
problem is reduced to the joint analysis of the following quantities:\vspace*{-1pt}
%
%e1.1 #&#
\begin{equation}
\label{u} U_n^{(r)}(t)=2^{-n\widetilde{\kappa}}\sum
_{k=0}^{\lfloor2^n t\rfloor
-1} H_r(Z_{(k+1)2^{-n}}-
Z_{k2^{-n}}),\qquad  n\in\mathbb{N}, t\geq0, r\in\N ^{*} .
\end{equation}
Here, $\widetilde{\kappa}>0$ is some constant depending a priori on
$r$, whereas $H_r$ denotes the $r$th Hermite polynomial ($H_1(x)=x$,
$H_2(x)=x^2-1$, etc.).
Due to the fact that one cannot separate $X$ from $Y$ inside $Z$ in the
definition
of $U_n^{(r)}$, working directly with (\ref{u}) seems to be a difficult
task (see also \cite{6}, Problem 5.1).
This is why, following an idea introduced by Khosnevisan and Lewis \cite
{5} in the study of the case $H=1/2$,
we will rather analyze $U_n^{(r)}$ by means of certain stopping times
for $Y$. The idea is quite simple: by stopping $Y$ as it crosses
certain levels, and by sampling $Z$ at these times, one can effectively
separate $X$ from $Y$.
To be more specific, let us introduce the following collection of
stopping times (with
respect to the natural filtration of $Y$), noted\vspace*{-1pt}
%
%e1.2 #&#
\begin{equation}\label{TN}
\mathscr{T}_n=\{T_{k,n}\dvt  k\geq0\},\qquad  n\geq0,
\end{equation}
which are in turn expressed in terms of the subsequent hitting
times of a dyadic grid cast on the real axis. More precisely, let
$\mathscr{D}_n= \{j2^{-n/2}\dvt  j\in\Z\}$, $n\geq0$, be the dyadic
partition (of $\R$) of order $n/2$. For every $n\geq0$, the
stopping times $T_{k,n}$, appearing in (\ref{TN}), are given by
the following recursive definition: $T_{0,n}= 0$, and\vspace*{-1pt}
\[
T_{k,n}= \inf \bigl\{s>T_{k-1,n}\dvt  Y(s)\in\mathscr{D}_n
\setminus\bigl\{Y(T_{k-1,n})\bigr\} \bigr\},\qquad  k\geq1.
\]
Note that the definition of $T_{k,n}$, and
therefore of $\mathscr{T}_n$, only involves the one-sided Brownian
motion $Y$, and that, for every $n\geq0$, the discrete stochastic
process\vspace*{-1pt}
\[
\mathscr{Y}_n=\bigl\{Y(T_{k,n})\dvt k\geq0\bigr\}
\]
defines a simple random
walk over $\mathscr{D}_n$. As shown in
\cite{5}, as $n$ tends to
infinity the collection $\{T_{k,n}\dvt  1\leq k \leq2^nt\}$ approximates the
common dyadic partition $\{k2^{-n}\dvt  1\leq k \leq2^nt\}$ of order $n$
of the time interval $[0,t]$ (see
\cite{5}, Lemma~2.2, for a precise statement).
Based on this fact, one can introduce the counterpart of (\ref{u})
based on $\mathscr{T}_n$, namely,\vspace*{-1pt}
%
%e1.3 #&#
\begin{equation}
\label{v} V_n^{(r)}(t)=2^{-n\widetilde{\kappa}}\sum
_{k=0}^{\lfloor2^n t\rfloor
-1} H_r(Z_{T_{k+1,n}}-
Z_{T_{k,n}}),\qquad  n\in\mathbb{N}, t\geq0, r\in\N^{*}.
\end{equation}

We are now in a position to state the main result of the present paper.

%th1.1 #&#
\begin{them}\label{main}
The following two f.d.d. convergences in law take place as $n\to\infty$
for any integer $N\geq1$.
\begin{longlist}[(2)]
\item[(1)] Assume that $H\leq\frac{1}2$. One has\vspace*{-1pt}
%
%e1.4 #&#
\begin{eqnarray}\label{conv1}
&&2^{-n/4} \Biggl\{ \sum_{k=0}^{\lfloor2^n t\rfloor-1}
H_{2r-1}\bigl(2^{nH/2}(Z_{T_{k+1,n}}- Z_{T_{k,n}})
\bigr)\dvt  1\leq r\leq N \Biggr\} _{t\geq0}
\nonumber
\\[-8pt]\\[-9pt]
&&\quad \stackrel{\mathrm{f.d.d.}} {\longrightarrow} \bigl\{ \sigma_{2r-1}
E^{(r)}_t \dvt  1\leq r\leq N \bigr\}_{t\geq0},\nonumber
\end{eqnarray}
where $\sigma_{2r-1}$ is some (explicit) constant and $E=(B^{(1)}\circ
Y,\ldots,B^{(N)}\circ Y)$, with $B=(B^{(1)},\ldots,\allowbreak B^{(N)})$ a
$N$-dimensional two-sided Brownian motion independent from $Y$.
\item[(2)] Assume that $H<\frac{3}4$. One has
%
%e1.5 #&#
\begin{eqnarray}\label{conv2}
&&2^{-3n/4} \Biggl\{ \sum_{k=0}^{\lfloor2^n t\rfloor-1}
H_{2r}\bigl(2^{nH/2}(Z_{T_{k+1,n}}- Z_{T_{k,n}})
\bigr)\dvt  1\leq r\leq N \Biggr\} _{t\geq0}
\nonumber
\\[-8pt]\\[-8pt]
&&\quad \stackrel{\mathrm{f.d.d.}} {\longrightarrow} \biggl\{ \sigma_{2r}
\int_{-\infty}^\infty L_t^x(Y)\,\mathrm{d}B^{(r)}_x
\dvt  1\leq r\leq N \biggr\}_{t\geq
0},\nonumber
\end{eqnarray}
where $\sigma_{2r}$ is some (explicit) constant, $B=(B^{(1)},\ldots
,B^{(N)})$ is a $N$-dimensional two-sided Brownian motion independent
from $Y$ and $L_t^x(Y)$ stands for the local time of $Y$ before time
$t$ at level~$x$.
\end{longlist}
\end{them}

The process $ \{ \int_{\R} L_t^x(Y)\,\mathrm{d}B_x^{(r)}\}_{t \geq0}$ appearing
in (\ref{conv2}) is nothing but the \textit{Brownian motion in Random
Scenery} introduced by Kesten and Spitzer (see \cite{9}).

As a corollary of this theorem, we deduce the fluctuations of the power
variation of $Z$.

%
%co1.2 #&#
\begin{cor}\label{cor-main}
The following two f.d.d. convergences in law take place as $n\to\infty
$ for any integer $N\geq1$.
\begin{longlist}[(2)]
\item[(1)] Assume that $H\leq\frac{1}2$. One has
%
%e1.6 #&#
\begin{eqnarray}\label
{cor-conv1}
&& \Biggl\{ 2^{\sklfrac{-n}{4}(1-(4r-2)H)}\sum_{k=0}^{\lfloor2^n t\rfloor
-1}
(Z_{T_{k+1,n}}- Z_{T_{k,n}})^{2r-1}\dvt  1\leq r\leq N \Biggr
\}_{t\geq
0}
\nonumber
\\[-8pt]\\[-8pt]
&&\quad \stackrel{\mathrm{f.d.d.}} {\longrightarrow} \Biggl\{ \sum
_{k=1}^r a_{r,k}\sigma_{2k-1}
E^{(k)}_t \dvt  1\leq r\leq N \Biggr\}_{t\geq0},\nonumber
\end{eqnarray}
where $a_{r,k}$ is some constant given by: ${a_{r,k} = \sum_{l=0}^{k-1}
\frac{(-1)^l (2(r+k-l-1))!}{l!(2(k-l)-1)!(r+k-l-1)!2^{r+k-1}}}$.
\item[(2)] Assume that $H<\frac{3}4$. One has
%
%e1.7 #&#
\begin{eqnarray}\label{cor-conv2}
&& \Biggl\{ 2^{\sklfrac{-3n}{4}(1-\sfrac{4rH}{3})}\sum_{k=0}^{\lfloor2^n
t\rfloor-1}
\bigl((Z_{T_{k+1,n}}- Z_{T_{k,n}})^{2r}- 2^{-nrH}b_{r,0}
\bigr)\dvt  1\leq r\leq N \Biggr\}_{t\geq0}\qquad
\nonumber\\[-8pt]\\[-8pt]
&&\quad \stackrel{\mathrm{f.d.d.}} {\longrightarrow} \Biggl\{ \sum
_{k=1}^r b_{r,k}\sigma_{2k}
\int_{-\infty}^\infty L_t^x(Y)\,\mathrm{d}B^{(k)}_x
\dvt  1\leq r\leq N \Biggr\}_{t\geq0},\nonumber
\end{eqnarray}
where $b_{r,k}$ is some constant given by: ${b_{r,k} = \sum_{l=0}^k
\frac{(-1)^l (2(r+k-l))!}{l!(2(k-l))!(r+k-l)!2^{r+k}}}$.
\end{longlist}
\end{cor}

Note that $b_{r,0}= E [(2^{nH/2}(Z_{T_{k+1,n}}-
Z_{T_{k,n}}))^{2r} ]= E[N^{2r}]$, with $N \sim\mathcal{N}(0,1)$.

In the particular case where $H=1/2$ (that is, when $Z$ is an iterated
Brownian motion) and $r=2,3,4$, we emphasize that Corollary~\ref
{cor-main} allows one to recovers Theorems 3.2, 4.4 and 4.5 from
Khoshnevisan and Lewis \cite{5}.

\begin{Remark*} To keep the length of this paper within bounds, I
defer to future analysis the technical investigation of the tightness
of the power variations of F.B.M.B.T given in the previous corollary.
\end{Remark*}

The organisation of the paper is as follows. In Section~\ref{sec2}, we provide
some needed preliminaries. Theorem~\ref{main} and Corollary~\ref
{cor-main} are then shown in Section~\ref{sec3}.

%s2 #&#
\section{Preliminaries}\label{sec2}
In this section, we collect several results that are useful for the
proof of Theorem~\ref{main}.

%s2.1 #&#
\subsection{An algebraic lemma and some local time estimates}\label{sec2.1}
For each
integer $n\geq0$, $k\in\Z$ and real number $t\geq0$, let $U_{j,n}(t)$ (resp.
$D_{j,n}(t)$) denote the number of \textit{upcrossings} (resp.
\textit{downcrossings}) of the interval
$[j2^{-n/2},(j+1)2^{-n/2}]$ within the first $\lfloor2^n
t\rfloor$ steps of the random walk $\{Y(T_{k,n})\}_{k\geq1}$, that is,
%
%e2.1 #&#
%e2.2 #&#
\begin{eqnarray}
\label{UPC}U_{j,n}(t)&=&\sharp \bigl\{k=0,\ldots,\bigl\lfloor2^nt\bigr
\rfloor-1 \dvt
\nonumber\\[-8pt]\\[-8pt]
&&\hphantom{\sharp \bigl\{}Y(T_{k,n}) = j2^{-n/2}\mbox{ and }Y(T_{k+1,n})=(j+1)2^{-n/2}
\bigr\};
\nonumber
\\
\label{DOWNC}D_{j,n}(t)&=&\sharp \bigl\{k=0,\ldots,\bigl\lfloor2^nt\bigr
\rfloor-1\dvt
\nonumber\\[-8pt]\\[-8pt]
&&\hphantom{\sharp \bigl\{}Y(T_{k,n}) = (j+1)2^{-n/2}\mbox{ and }Y(T_{k+1,n})=j2^{-n/2}
\bigr\}.
\nonumber
\end{eqnarray}
The following lemma will play a crucial role in our
study of the asymptotic behavior of $V_n^{(r)}$.
Its main feature is to
separate $X$ from $Y$, thus providing a representation of
$V_n^{(r)}$ which is amenable to analysis.

%le2.1 #&#
\begin{lemme}[(See \cite{5}, Lemma~2.4)]\label{lm-kl}  Fix $t\geq0$ and
$r\in\N^{*}$.
Then
%
%e2.3 #&#
\begin{eqnarray}\label{amenable}
V_n^{(r)}(t)= 2^{-n\widetilde{\kappa}}\sum
_{j\in\Z} H_r \bigl( 2^{nH/2}
(X_{(j+1)2^{-\sfrac{n}{2}}}- X_{j2^{-\sfrac{n}{2}}} ) \bigr) \bigl(U_{j,n}(t)+(-1)^r
D_{j,n}(t) \bigr).\qquad
\end{eqnarray}
\end{lemme}

Also, in order to prove the second point of Theorem~\ref{main} we will
need estimates on the local time of $Y$
taken from \cite{5}, that we collect in the following statement.

%
%pr2.2 #&#
\begin{prop}\label{properties-kl}
\begin{longlist}[3.]
\item[1.] For every $x\in\R$, $p \in\N^{*}$ and $t > 0$, we have
\[
E \bigl[\bigl( L_t^{x}(Y)\bigr)^p \bigr]
\leq2 E \bigl[\bigl(L_1^0(Y)\bigr)^p \bigr]
t^{p/2} \exp \biggl(-\frac{x^2}{2t} \biggr).
\]
\item[2.] There exists a positive constant $\mu$ such that, for every
$a,b\in\R$ with $ab\geq0$
and $t > 0$,
\[
E \bigl[\bigl|L_t^b(Y)-L_t^a(Y)\bigr|^2
\bigr]^{1/2}\leq\mu\sqrt{|b-a|} t^{1/4} \exp \biggl(-
\frac{a^2}{4t} \biggr).
\]
\item[3.] There exists a positive random variable $K\in L^8$ such that, for
every $j\in\Z$, every $n\geq0$
and every $t > 0$, one has that
\[
\bigl|\mathcal{L}_{j,n}(t)-L_t^{j2^{-n/2}}(Y)\bigr|
\leq2Kn2^{-n/4}\sqrt {L_t^{j2^{-n/2}}(Y)},
\]
where $\mathcal{L}_{j,n}(t)=2^{-n/2}(U_{j,n}(t)+D_{j,n}(t))$.
\end{longlist}
\end{prop}

%s2.2 #&#
\subsection{Breuer--Major Theorem}\label{sec2.2}

Let $\lbrace G_{k}\rbrace_{k\geq1}$ be a centered stationary Gaussian
sequence. In this Gaussian context, stationary just means that there
exist $\rho\dvtx  \mathbb{Z}\rightarrow\mathbb{R}$ such that $E[G_{k}G_{l}]
= \rho(k-l)$, $k,l\geq1$. Assume further that $\rho(0)=1$, that is,
each $G_{k}$ is $\mathscr{N}(0,1)$ distributed. Let $\varphi\dvtx  \mathbb
{R}\rightarrow\mathbb{R}$ be a measurable function satisfying
%
%e2.4 #&#
\begin{eqnarray}
\label{1} E\bigl[\varphi^{2}(G_{1})\bigr]=
\frac{1}{\sqrt{2\pi}}\int_{\mathbb{R}}\varphi ^{2}(x)\mathrm{e}^{-x^{2}/2}\,\mathrm{d}x
< +\infty.
\end{eqnarray}
The function $\varphi$ may be expanded in $L^{2}(\mathbb{R}, \frac
{\mathrm{e}^{-x^{2}/2}}{\sqrt{2 \pi}}\,\mathrm{d}x)$ (in a unique way) in terms of Hermite
polynomials as follows:
%
%e2.5 #&#
\begin{eqnarray}
\label{2} \varphi(x)=\sum_{q=0}^{+\infty}a_{q}H_{q}(x).
\end{eqnarray}
Let $d\geq0$ be the \textit{Hermite rank} of $\varphi$, that is, the
first integer $q\geq0$ such that $a_{q}\neq0$ in (\ref{2}). We then
have the celebrated Breuer--Major Theorem (see \cite{11}, see also \cite
{8} for a modern proof).

%
%th2.3 #&#
\begin{theorem}[(Breuer--Major)]\label{th:bm}
Let $\lbrace G_{k} \rbrace_{k\geq1}$ (with covariance $\rho$) and
$\varphi\dvtx \mathbb{R}\rightarrow\mathbb{R}$ (with Hermite index $d$) be
as above. Assume further that $\sum_{k\in\mathbb{Z}}|\rho(k)|^{d} <
+\infty$. Then, as $n\rightarrow+\infty$,
%
%e2.6 #&#
\begin{eqnarray}
\label{3} 2^{-n/2} \Biggl\{\sum_{k=1}^{\lfloor2^nt\rfloor}
\bigl(\varphi(G_{k}) - E\bigl[ \varphi(G_{k})\bigr]
\bigr) \Biggr\}_{t\geq0}\stackrel{\mathrm{f.d.d.}} {\to } \{\sigma
B_t\}_{t\geq0},
\end{eqnarray}
with $B$ a standard Brownian motion and $\sigma>0$ given by
%
%e2.7 #&#
\begin{eqnarray}
\label{4} \sigma^{2}=\sum_{q=d}^{+\infty}q!a_{q}^{2}
\sum_{k\in\mathbb{Z}}\rho (k)^{q} \in[0,+\infty[.
\end{eqnarray}
\end{theorem}

%s2.3 #&#
\subsection{Peccati--Tudor Theorem}\label{sec2.3}
In a seminal paper of 2005, Nualart and Peccati \cite{12} discovered a
surprising central limit theorem (called the \textit{Fourth Moment
Theorem} nowadays) for sequences of multiple stochastic integrals of a
fixed order: in this context, convergence in distribution to the
standard normal law is actually equivalent to convergence of just the
fourth moment. Shortly afterwards, Peccati and Tudor gave a
multidimensional version of this characterization, making use of tools
belonging to the Malliavin calculus. Since we will rely on this result
in the present paper, let us give more details.

Let $d\geq2$ and $q_1,\ldots,q_d\geq1$ be some fixed integers.
Consider a sequence of random vectors $F_n=(F_{1,n},\ldots,F_{d,n})$ of
the following form.
Each $F_{i,n}$ can be written as
\[
F_{i,n}=\sum_{j=0}^{N_{n}}
a_{j,n} H_{q_i}(Y_{j,n}),
\]
where $N_{n}$ is an integer, $a_{j,n}$ are real numbers and $\lbrace
Y_{j,n}\rbrace_{j\geq0}$ is a centered stationary Gaussian family with
unit variance.
We then have the following result, shown in \cite{10}.

%
%th2.4 #&#
\begin{theorem}[(Peccati--Tudor)]\label{th:pec}
Let $(F_n)$ be a sequence as above. Let $C\in\mathscr{M}_{d}(\mathbb
{R})$ be a symmetric and positive matrix, and let N be a centered
Gaussian vector with covariance $C$. Assume that
%
%e2.8 #&#
\begin{eqnarray}
\label{5} \lim_{n\rightarrow+\infty}E[F_{i,n}F_{j,n}]=C(i,j),
\qquad 1\leq i,j \leq d.
\end{eqnarray}
Then, as $n\rightarrow+\infty$, the following two conditions are
equivalent:
\begin{enumerate}[(b)]
\item[(a)] $F_{n}$ converges in law to $N$;

\item[(b)] for every $1\leq i \leq d$, $F_{i,n}$ converges in law to $\mathscr
{N}(0,C(i,i))$.
\end{enumerate}
\end{theorem}

%s3 #&#
\section{Proof of Theorem \texorpdfstring{\protect\ref{main}}{1.1}}\label{sec3}

%s3.1 #&#
\subsection{Proof of \texorpdfstring{(\protect\ref{conv1})}{(1.4)}}\label{sec3.1}

Recall the definition (\ref{v}) of $V_n^{(r)}(t)$ and let us fix
$\widetilde{\kappa}=1/4$.
First of all, let us apply Lemma~\ref{lm-kl}. Because $2r-1$ is an odd
number, we obtain that
%
%e3.1 #&#
\begin{equation}
V_n^{(2r-1)}(t)=2^{-n/4} \sum
_{j\in\mathbb{Z}} H_{2r-1}\bigl(2^{nH/2}(X_{(j+1)2^{-n/2}}-
X_{j2^{-n/2}})\bigr) \bigl(U_{j,n}(t)-D_{j,n}(t)\bigr).
\end{equation}
Now, let us observe (see also \cite{5}, Lemma~2.5) that
\[
U_{j,n}(t) - D_{j,n}(t)= \lleft\{ %
\begin{array}
{l@{\qquad}l} \mathbf{1}\bigl(0\leq j< j^*(n,t)\bigr) & \mbox{if }j^*(n,t) > 0,
\\
0 & \mbox{if }j^{*} = 0,
\\
-\mathbf{1}\bigl(j^*(n,t)\leq j<0\bigr) & \mbox{if }j^*(n,t)< 0, \end{array}
\rright.
\]
where
\[
j^*(n,t)=2^{n/2}Y_{T_{\lfloor2^n t\rfloor,n}}.
\]
As a consequence,
\begin{eqnarray*}
V_n^{(2r-1)}(t)=\lleft\{ %
\begin{array} {l@{\qquad}l}
2^{-n/4} \displaystyle \sum_{j=1}^{j^*(n,t)}
H_{2r-1}\bigl(2^{nH/2}\bigl(X^+_{j2^{-n/2}}-
X^+_{(j-1)2^{-n/2}}\bigr)\bigr) & \mbox{if }j^*(n,t) > 0,
\\\noalign{\vspace*{2pt}}
0 &\mbox{if }j^{*} = 0,
\\\noalign{\vspace*{2pt}}
2^{-n/4} \displaystyle \sum_{j=1}^{|j^*(n,t)|}
H_{2r-1}\bigl(2^{nH/2}\bigl(X^-_{j2^{-n/2}}-
X^-_{(j-1)2^{-n/2}}\bigr)\bigr) & \mbox{if }j^*(n,t) < 0, \end{array}
\rright.
\end{eqnarray*}
where $X^+_t=X_t$ for $t\geq0$ and $X^-_{-t}=X_{t}$ for $t<0$.
Our analysis of $V_n^{(2r-1)}$ will become easier if one
introduces the following sequence of processes $W_{\pm,n}^{(2r-1)}$, in
which we have replaced
$\sum_{j=1}^{\pm j^*(n,t)}$ by $\sum_{j=1}^{\lfloor2^{n/2}t\rfloor}$, namely:
\begin{eqnarray*}
W_{+,n}^{(2r-1)}(t)&=& 2^{-n/4} \sum
_{j=1}^{\lfloor2^{n/2}t\rfloor} H_{2r-1}\bigl(2^{nH/2}
\bigl(X^+_{j2^{-n/2}}- X^+_{(j-1)2^{-n/2}}\bigr)\bigr),\qquad  t \geq0,
\\
W_{-,n}^{(2r-1)}(t)&=& 2^{-n/4} \sum
_{j=1}^{\lfloor2^{n/2}t\rfloor} H_{2r-1}\bigl(2^{nH/2}
\bigl(X^-_{j2^{-n/2}}- X^-_{(j-1)2^{-n/2}}\bigr)\bigr) ,\qquad  t \geq0,
\\
W_{n}^{(2r-1)}(t)&=&\lleft\{ %
\begin{array} {l@{\qquad}l}
W_{+,n}^{(2r-1)}(t) &\mbox{if }t \geq0,
\\\noalign{\vspace*{4pt}}
W_{-,n}^{(2r-1)}(-t) &\mbox{if }t < 0. \end{array} %
\rright.
\end{eqnarray*}
It is clear, using the self-similarity property of $X$, that the f.d.d.
convergence in law of the vector $(W_{+,n}^{(2r-1)},W_{-,n}^{(2r-1)},
1\leq r\leq N)$ is equivalent to the f.d.d. convergence in law of the
vector $(\overline{W}_{+,n}^{(2r-1)},\overline{W}_{-,n}^{(2r-1)}, 1\leq
r\leq N)$ defined as:
\begin{eqnarray*}
\overline{W}_{+,n}^{(2r-1)}(t)&=& 2^{-n/4} \sum
_{j=1}^{\lfloor2^{n/2}t\rfloor} H_{2r-1}\bigl(X^+_{j}-
X^+_{j-1}\bigr), \qquad t \geq0,
\\
\overline{W}_{-,n}^{(2r-1)}(t)&=& 2^{-n/4} \sum
_{j=1}^{\lfloor2^{n/2}t\rfloor} H_{2r-1}\bigl(X^-_{j}-
X^-_{j-1}\bigr) , \qquad t \geq0.
\end{eqnarray*}
Let $G_j=X^+_j-X^+_{j-1}$. The family $\{G_j\}$ is Gaussian,
stationary, centered, with variance 1;
moreover its covariance $\rho$ is given by
%
%e3.2 #&#
\begin{equation}
\label{rho} \rho(k)=E[G_jG_{j+k}]=\tfrac{1}{2}
\bigl(|k+1|^{2H}+|k-1|^{2H}-2|k|^{2H} \bigr),
\end{equation}
so that $\sum|\rho(k)|<\infty$ because $H\leq\frac{1}2$.
Hence, Breuer--Major Theorem \ref{th:bm} applies and yields that, as
$n\to\infty$ and for any fixed $r$,
\[
\bigl\{\overline{W}_{+,n}^{(2r-1)}(t)\dvt  t\geq0\bigr\}\stackrel{
\mathrm {f.d.d.}} {\longrightarrow}\sigma_{2r-1}\bigl\{B^{+,r}(t)\dvt
t\geq0\bigr\},
\]
with $B^{+,r}$ a standard Brownian motion and $\sigma_{2r-1}=\sqrt {(2r-1)!\sum_{a\in\Z}\rho(a)^{2r-1}}$. Note that
$\sum_{a\in\Z}|\rho(a)|^{2r-1} < \infty$ if and only if $H < 1-
1/(2(2r-1))$, which is satisfied for all $r \geq1 $ since we have
supposed that $H \leq1/2$ (the case $H=1/2$ may be treated separately).
Similarly,
\[
\bigl\{\overline{W}_{-,n}^{(2r-1)}(t)\dvt  t\geq0\bigr\}\stackrel{
\mathrm {f.d.d.}} {\longrightarrow}\sigma_{2r-1}\bigl\{B^{-,r}(t)\dvt
t\geq0\bigr\},
\]
with $B^{-,r}$ a standard Brownian motion and $\sigma_{2r-1}$ as above.
In order to deduce the joint convergence in law of $(\overline
{W}_{+,n}^{(2r-1)},\overline{W}_{-,n}^{(2r-1)}, 1\leq r\leq N)$, from
Peccati--Tudor Theorem \ref{th:pec} and taking into account that
$E[\overline{W}_{\pm,n}^{(2r-1)}(t)\overline{W}_{\pm,n}^{(2l-1)}(t)]=0$
for $l\neq r$ (since Hermite polynomials of different orders are
orthogonal), it remains to check that, for any integer $r$ and any real
numbers $t , s \geq0$,
%
%e3.3 #&#
\begin{equation}
\label{limit-to-check} \lim_{n\to\infty} E\bigl[\overline{W}_{+,n}^{(2r-1)}(t)
\overline {W}_{-,n}^{(2r-1)}(s)\bigr]=0.
\end{equation}
Let us do it. One can write,
\begin{eqnarray*}
&&E\bigl[\overline{W}_{+,n}^{(2r-1)}(t)\overline{W}_{-,n}^{(2r-1)}(s)
\bigr]
\\
&&\quad = 2^{-n/2}\sum_{k=1}^{\lfloor2^{n/2}t\rfloor}\sum
_{l=1}^{\lfloor
2^{n/2}s\rfloor} E\bigl[H_{2r-1}
\bigl(X^+_{k}- X^+_{k-1}\bigr)H_{2r-1}
\bigl(X^-_{l}- X^-_{l-1}\bigr)\bigr]
\\
&&\quad =(2r-1)! 2^{-n/2}\sum_{k=1}^{\lfloor2^{n/2}t\rfloor}
\sum_{l=1}^{\lfloor2^{n/2}s\rfloor} \bigl(E\bigl[
\bigl(X^+_{k}- X^+_{k-1}\bigr) \bigl(X^-_{l}-
X^-_{l-1}\bigr)\bigr] \bigr)^{2r-1}
\\
&&\quad =(2r-1)! 2^{-n/2}\sum_{k=1}^{\lfloor2^{n/2}t\rfloor}
\sum_{l=1}^{\lfloor2^{n/2}s\rfloor} \bigl(E\bigl[(X_{k}-
X_{k-1}) (X_{-l}- X_{-l+1})\bigr]
\bigr)^{2r-1}
\\
&&\quad =(2r-1)! 2^{-n/2}\sum_{k=1}^{\lfloor2^{n/2}t\rfloor}
\sum_{l=1}^{\lfloor2^{n/2}s\rfloor} \biggl( \frac{1}{2}
\bigl[ 2|k+l-1|^{2H} - |k+l|^{2H} - |k+l-2|^{2H}
\bigr] \biggr)^{2r-1}.
\end{eqnarray*}
Setting $a=k+l$, we deduce that
\begin{eqnarray*}
&&\hspace*{-2pt}E\bigl[\overline{W}_{+,n}^{(2r-1)}(t)\overline{W}_{-,n}^{(2r-1)}(s)
\bigr]
\\
&&\hspace*{-2pt}\quad = 2^{-(2r-1)}(2r-1)! 2^{-n/2}\sum_{k=1}^{\lfloor2^{n/2}t\rfloor}
\sum_{a=k+1}^{k+\lfloor2^{n/2}s\rfloor} \bigl(2|a-1|^{2H}-
|a|^{2H} -|a-2|^{2H} \bigr)^{2r-1}
\\
&&\hspace*{-2pt}\quad =2^{-(2r-1)}(2r-1)! 2^{-n/2}\sum_{a=2}^{\lfloor2^{n/2}t\rfloor+
\lfloor2^{n/2}s\rfloor}
\sum_{k=(a-\lfloor2^{n/2}s\rfloor)\vee
1}^{(a-1)\wedge\lfloor2^{n/2}t\rfloor} \bigl(2|a-1|^{2H}-
|a|^{2H} -|a-2|^{2H} \bigr)^{2r-1}
\\
&&\hspace*{-2pt}\quad =2^{-(2r-1)}(2r-1)! 2^{-n/2}\sum_{a\in\N}f_n(a),
\end{eqnarray*}
where
\begin{eqnarray*}
f_n(a)&:=& \bigl(2|a-1|^{2H}- |a|^{2H}
-|a-2|^{2H} \bigr)^{2r-1} \bigl( (a-1)\wedge\bigl
\lfloor2^{n/2}t\bigr\rfloor- \bigl(a - \bigl\lfloor2^{n/2}s\bigr
\rfloor\bigr)\vee 1 +1 \bigr)
\\
&&{}\times\1_{\{2 \leq a \leq\lfloor2^{n/2}t\rfloor+ \lfloor
2^{n/2}s\rfloor\}} .
\end{eqnarray*}
For any $a\in\{2,\ldots,\lfloor2^{n/2}t\rfloor+ \lfloor
2^{n/2}s\rfloor\}$,
observe that $\lfloor2^{n/2}t\rfloor+ \lfloor2^{n/2}s\rfloor\geq1$.
Also, we have
\begin{eqnarray*}
 2^{-\sfrac{n}{2}} \bigl| (a-1)\wedge\bigl\lfloor2^{n/2}t\bigr\rfloor \bigr| \1
_{\{2 \leq a \leq\lfloor2^{n/2}t\rfloor+ \lfloor2^{n/2}s\rfloor\}} \leq2^{-\sfrac{n}{2}}\bigl\lfloor2^{n/2}t\bigr\rfloor
\leq t,
\end{eqnarray*}
as well as
\begin{eqnarray*}
&& 2^{-\sfrac{n}{2}}\bigl |\bigl(a - \bigl\lfloor2^{n/2}s\bigr\rfloor\bigr)
\vee1 \bigr|\1_{\{
2 \leq a \leq\lfloor2^{n/2}t\rfloor+ \lfloor2^{n/2}s\rfloor\}}
\\
&&\quad \leq 2^{-\sfrac{n}{2}} \bigl( |a|+ \bigl\lfloor2^{n/2}s\bigr\rfloor+ 1
\bigr)\1_{\{2
\leq a \leq\lfloor2^{n/2}t\rfloor+ \lfloor2^{n/2}s\rfloor\}}
\\
&&\quad \leq 2^{-\sfrac{n}{2}}\bigl(\bigl\lfloor2^{n/2}t\bigr\rfloor+ 2\bigl
\lfloor 2^{n/2}s\bigr\rfloor\bigr) + 2^{-\sfrac{n}{2}}\1_{\{2 \leq a \leq\lfloor
2^{n/2}t\rfloor+ \lfloor2^{n/2}s\rfloor\}}
\\
&&\quad  \leq 2^{-\sfrac{n}{2}}\bigl(\bigl\lfloor2^{n/2}t\bigr\rfloor+ 2\bigl
\lfloor 2^{n/2}s\bigr\rfloor\bigr) + 2^{-\sfrac{n}{2}}\bigl(\bigl
\lfloor2^{n/2}t\bigr\rfloor+ \bigl\lfloor2^{n/2}s\bigr\rfloor
\bigr)\leq2t+3s,
\end{eqnarray*}
and
\begin{eqnarray*}
 2^{-\sfrac{n}{2}}\1_{\{2 \leq a \leq\lfloor2^{n/2}t\rfloor+
\lfloor2^{n/2}s\rfloor\}} \leq2^{-\sfrac{n}{2}} \bigl( \bigl\lfloor
2^{n/2}t\bigr\rfloor+ \bigl\lfloor2^{n/2}s\bigr\rfloor\bigr)
\leq t+s.
\end{eqnarray*}
Plugging all these inequalities together leads to
\[
2^{-n/2}\bigl|f_n(a)\bigr| \leq(4t+4s) \bigl|2|a-1|^{2H}-
|a|^{2H} -|a-2|^{2H} \bigr|^{2r-1}
\]
for all $n$,
with $\sum_{a\in\mathbb{N}}  |2|a-1|^{2H}- |a|^{2H} -|a-2|^{2H}
|^{2r-1}<\infty$ (recall that $H \leq1/2$).
Moreover, ${ 2^{-n/2}f_n(a)\displaystyle\mathop{\longrightarrow}_{n
\rightarrow\infty} 0}$ for any fixed $a$ because
\[
\frac{(a-1)}{2^{n/2}}\wedge\frac{\lfloor2^{n/2}t\rfloor}{2^{n/2}} - \frac{(a - \lfloor2^{n/2}s\rfloor)}{2^{n/2}}
\vee2^{-n/2} + 2^{-n/2} \displaystyle\mathop{\longrightarrow}_{n \to\infty}
0 \wedge t - (-s)\vee0 = 0
\]
since $t,s \geq0$.
Hence, the dominated convergence theorem applies and yields
\[
2^{-n/2}\sum_{a\in\N}f_n(a)
\displaystyle\mathop{\longrightarrow}_{n
\rightarrow\infty} 0,
\]
that is, (\ref{limit-to-check}) holds true. As we said, using
Peccati--Tudor Theorem \ref{th:pec} one thus obtains that
%
%e3.4 #&#
\begin{equation}
\label{cv1} \bigl(W_{+,n}^{(2r-1)},W_{-,n}^{(2r-1)},
1\leq r\leq N\bigr) \stackrel{\mathrm{f.d.d.}} {\longrightarrow} \bigl(
\sigma_{2r-1}B^{+,r},\sigma_{2r-1}B^{-,r}, 1\leq
r\leq N\bigr),
\end{equation}
with $(B^{+,r},B^{-,r}, 1\leq r\leq N)$ a $2N$-dimensional standard
Brownian motion. As a consequence, we have
%
%e3.5 #&#
\begin{equation}
\label{cv2} \bigl(W_{n}^{(2r-1)}(t), 1\leq r\leq N
\bigr)_{t \in\R} \stackrel{\mathrm{f.d.d.}} {\longrightarrow} \bigl(
\sigma_{2r-1}B^{(r)}(t), 1\leq r\leq N\bigr)_{t \in\R},
\end{equation}
with $(B^{(r)}, 1\leq r\leq N)$ a $N$-dimensional two-sided Brownian motion.

On the other hand, let us prove for any $r\in\N^{*}$ the existence of $C_r>0$
such that, for any $n$ and any $s,t\in\R$,
%
%e3.6 #&#
\begin{equation}
\label{bound} E\bigl[\bigl(W_n^{(2r-1)}(t) -
W_n^{(2r-1)}(s)\bigr)^2\bigr]\leq8C_r
\bigl(2^{-n/2} + |t-s|\bigr).
\end{equation}

To do so, we distinguish three cases, according to the sign of $s, t
\in\R$ (and reducing the problem by symmetry):
\begin{longlist}[(3)]
\item[(1)] if $0 \leq s \leq t$:
\begin{eqnarray*}
&&E\bigl[\bigl(W_n^{(2r-1)}(t) - W_n^{(2r-1)}(s)
\bigr)^2\bigr]
\\
&&\quad = E\bigl[\bigl(\overline{W}_{+,n}^{(2r-1)}(t)-\overline
{W}_{+,n}^{(2r-1)}(s)\bigr)^2\bigr]=
\Biggl|2^{-n/2}E\Biggl[ \Biggl( \sum_{j=\lfloor2^{n/2}
s\rfloor+1}^{\lfloor2^{n/2} t\rfloor}
H_{2r-1}\bigl(X_j^{+}-X_{j-1}^{+}
\bigr) \Biggr)^2\Biggr] \Biggr|
\\
&&\quad = \Biggl|(2r-1)! 2^{-n/2}\sum_{j,k=\lfloor2^n s\rfloor+1}^{\lfloor2^n
t\rfloor}
\bigl(\rho(j-k)\bigr)^{2r-1} \Biggr|
\\
&&\quad = \Biggl|(2r-1)! 2^{-n/2}\sum_{j=\lfloor2^{n/2} s\rfloor+1}^{\lfloor
2^{n/2} t\rfloor}
\sum_{a=j-\lfloor2^{n/2} t\rfloor}^{j-\lfloor
2^{n/2} s\rfloor-1}\bigl(\rho(a)
\bigr)^{2r-1} \Biggr|
\\
&&\quad \leq(2r-1)! 2^{-n/2}\\
&&\qquad {}\times\sum_{a=\lfloor2^{n/2} s\rfloor-\lfloor2^{n/2}
t\rfloor+1}^{\lfloor2^{n/2} t\rfloor-\lfloor2^{n/2} s\rfloor
-1}\bigl|
\rho(a)\bigr|^{2r-1}
\bigl |\bigl(a + \bigl\lfloor2^n t\bigr\rfloor\bigr)\wedge\bigl
\lfloor2^n t\bigr\rfloor- \bigl(a+ \bigl\lfloor2^n s\bigr
\rfloor\bigr)\vee\bigl(\bigl\lfloor2^n s\bigr\rfloor\bigr)\bigr |
\\
&&\quad  \leq(2r-1)! 2^{-n/2}\sum_{a \in\Z}\bigl|
\rho(a)\bigr|^{2r-1} \bigl|\bigl\lfloor 2^{n/2} t\bigr\rfloor- \bigl
\lfloor2^{n/2} s \bigr\rfloor \bigr|= C_r 2^{-n/2} \bigl|\bigl
\lfloor2^{n/2} t\bigr\rfloor- \bigl\lfloor2^{n/2} s \bigr\rfloor
\bigr|
\\
&&\quad  \leq C_r \bigl(\bigl|2^{-n/2}\bigl\lfloor2^{n/2} t
\bigr\rfloor- t\bigr| + |t-s| + \bigl|2^{-n/2}\bigl\lfloor2^{n/2} s \bigr
\rfloor-s\bigr|\bigr) \leq C_r\bigl(22^{-n/2} + |t-s|\bigr)
\\
&&\quad  \leq2C_r\bigl(2^{-n/2} +|t-s|\bigr),
\end{eqnarray*}
with ${C_r= (2r-1)!\sum_{a \in\Z}|\rho(a)|^{2r-1}}< \infty$, hence
(\ref{bound}) holds true.
\item[(2)] If $s \leq t \leq0$: by the same argument as above
\begin{eqnarray*}
&&E\bigl[\bigl(W_n^{(2r-1)}(t) - W_n^{(2r-1)}(s)
\bigr)^2\bigr]
\\
&&\quad =E\bigl[\bigl(\overline{W}_{-,n}^{(2r-1)}(-s)- \overline
{W}_{-,n}^{(2r-1)}(-t)\bigr)^2\bigr]\leq
C_r\bigl(2^{-n/2} \bigl|\bigl\lfloor2^{n/2} |t|\bigr
\rfloor- \bigl\lfloor2^{n/2} |s|\bigr\rfloor \bigr|\bigr)
\\
&&\quad  \leq C_r \bigl( \bigl|2^{-n/2}\bigl\lfloor2^{n/2} |t|
\bigr\rfloor- |t| \bigr| + \bigl||t|-|s| \bigr| + \bigl|2^{-n/2}\bigl\lfloor2^{n/2}
|s| \bigr\rfloor-|s| \bigr|\bigr)
\\
&&\quad  \leq C_r\bigl(22^{-n/2} + \bigl||t|-|s| \bigr|\bigr)
\\
&&\quad  \leq2C_r\bigl(2^{-n/2} + \bigl||t|-|s| \bigr|\bigr)
\leq2C_r\bigl(2^{-n/2} +|t-s|\bigr),
\end{eqnarray*}
so that (\ref{bound}) holds true as well.
\item[(3)] If $ s < 0 < t$: using the two previous inequality (point
(1) and point (2)) one has
\begin{eqnarray*}
E\bigl[\bigl(W_n^{(2r-1)}(t) - W_n^{(2r-1)}(s)
\bigr)^2\bigr]&\leq &2E\bigl[\bigl(W_{n}^{(2r-1)}(t)-W_{n}^{(2r-1)}(0)
\bigr)^2\bigr]
\\
&&{} + 2E\bigl[\bigl(W_{n}^{(2r-1)}(s)-W_{n}^{(2r-1)}(0)
\bigr)^2\bigr]
\\
& \leq&4C_r\bigl(2^{-n/2} +t\bigr) + 4C_r
\bigl(2^{-n/2} + |s|\bigr)
\\
&=& 8C_r2^{-n/2} + 4C_r\bigl(t + |s|\bigr)
\\
&=& 8C_r2^{-n/2}+4C_r|t -s|
\\
&\leq&8C_r\bigl(2^{-n/2} +|t-s|\bigr).
\end{eqnarray*}
\end{longlist}
This proves (\ref{bound}).

Now, let us go back to $V_n^{(2r-1)}$. Observe that
\[
V_n^{(2r-1)}(t) = W_n^{(2r-1)}(Y_{T_{\lfloor2^n t\rfloor,n}})
\]
and recall from \cite{5}, Lemma~2.3, that $E[|Y_{T_{\lfloor2^n t\rfloor
,n}}-Y(t)|]\to0$ as $n\to\infty$ for
any $t > 0$.
We deduce, combining these two latter facts with (\ref{bound}), that
\[
V_n^{(2r-1)}(t) - W_n^{(2r-1)}(Y_{t})
\stackrel{L^2} {\to} 0 \qquad \mbox{as }n\to\infty.
\]
But $Y$ is independent from $W_n^{(2r-1)}$ so, from (\ref{cv2}), it
comes that
%
%e3.7 #&#
\begin{equation}
\label{cv3} \bigl(W_{n}^{(2r-1)}\circ Y, 1\leq r\leq N\bigr)
\stackrel{\mathrm{f.d.d.}} {\longrightarrow} \bigl(\sigma_{2r-1}B^{(r)}
\circ Y, 1\leq r\leq N\bigr),
\end{equation}
from which the desired conclusion (\ref{conv1}) follows.

%s3.2 #&#
\subsection{Proof of \texorpdfstring{(\protect\ref{conv2})}{(1.5)}}\label{sec3.2}

Recall the definition (\ref{v}) and let us fix $\widetilde{\kappa}=3/4$.
First of all, let us apply Lemma~\ref{lm-kl}. Because $2r$ is an even
number, we obtain that
%
%e3.8 #&#
\begin{equation}
V_n^{(2r)}(t)=2^{-3n/4} \sum
_{j\in\mathbb{Z}} H_{2r}\bigl(2^{nH/2}(X_{(j+1)2^{-n/2}}-
X_{j2^{-n/2}})\bigr) \bigl(U_{j,n}(t)+D_{j,n}(t)\bigr).
\end{equation}
Set
$
\mathcal{L}_{j,n}(t)=2^{-n/2}(U_{j,n}(t) + D_{j,n}(t))$,
so that
\begin{eqnarray*}
V_n^{(2r)}(t)=2^{-n/4} \sum
_{j\in\mathbb{Z}} H_{2r}\bigl(2^{nH/2}(X_{(j+1)2^{-n/2}}-
X_{j2^{-n/2}})\bigr)\mathcal{L}_{j,n}(t).
\end{eqnarray*}
At this stage, to simplify the exposition, let us introduce the
short-hand notation
\[
X^{(n)}_j=2^{nH/2}X_{j2^{-n/2}}.
\]

Fix $t\geq0$. In order to study the convergence in law of
$V_n^{(2r)}(t)$ as $n$ tends to infinity, we shall consider
(separately) the cases when $n$ is even and when $n$ is odd.

When $n$ is even, for any even integers $n\geq m\geq0$ and any
integer $p\geq0$, by following Nourdin and Peccati (see \cite{7}) one
can decompose $V_n^{(2r)}(t)$ as
\[
V_n^{(2r)}(t) = A_{m,n,p}^{(2r)}(t) +
B_{m,n,p}^{(2r)}(t) + C_{m,n,p}^{(2r)}(t) +
D_{m,n,p}^{(2r)}(t) ,
\]
where
\begin{eqnarray*}
A_{m,n,p}^{(2r)}(t) &=&2^{-n/4} \sum
_{-p 2^{m/2} +1\leq j\leq p 2^{m/2}} \sum_{i=(j-1)2^{\vfrac
{n-m}{2}}}^{j2^{\vfrac{n-m}{2}}-1}
H_{2r}\bigl(X^{(n)}_{i+1}-X^{(n)}_i
\bigr) \bigl(\mathcal{L}_{i,n}(t) - L_t^{i2^{-n/2}}(Y)
\bigr),
\\
B_{m,n,p}^{(2r)}(t) &=&2^{-n/4} \sum
_{-p 2^{m/2}+1\leq j\leq p 2^{m/2}} \sum_{i=(j-1)2^{\vfrac
{n-m}{2}}}^{j2^{\vfrac{n-m}{2}}-1}
H_{2r}\bigl(X^{(n)}_{i+1}-X^{(n)}_i
\bigr)
\\
&&\hphantom{2^{-n/4} \sum
_{-p 2^{m/2}+1\leq j\leq p 2^{m/2}} \sum_{i=(j-1)2^{\vfrac
{n-m}{2}}}^{j2^{\vfrac{n-m}{2}}-1}}{} \times\bigl(L_t^{i2^{-n/2}}(Y)-L_t^{j2^{-m/2}}(Y)
\bigr),
\\
C_{m,n,p}^{(2r)}(t) &=&2^{-n/4} \sum
_{-p 2^{m/2}+1\leq j\leq p 2^{m/2}}L_t^{j2^{-m/2}}(Y) \sum
_{i=(j-1)2^{\vfrac{n-m}{2}}}^{j2^{\vfrac{n-m}{2}}-1} H_{2r}\bigl(X^{(n)}_{i+1}-X^{(n)}_i
\bigr),
\\
D_{m,n,p}^{(2r)}(t) &=&2^{-n/4} \sum
_{i\geq p 2^{n/2}} H_{2r}\bigl(X^{(n)}_{i+1}-X^{(n)}_i
\bigr)\mathcal {L}_{i,n}(t) + \sum_{i < -p2^{n/2}}
H_{2r}\bigl(X^{(n)}_{i+1}-X^{(n)}_i
\bigr)\mathcal{L}_{i,n}(t).
\end{eqnarray*}
We can see that since we have taken even integers $n\geq m\geq0$ then
$2^{m/2}$, $2^{\vfrac{n-m}{2}}$ and $2^{n/2}$ are integers as well. This
justifies the validity of the previous decomposition.

When $n$ is odd, for any odd integers $n \geq m \geq0$ we can work
with the same decomposition for $V_n^{(2r)}(t)$. The only difference is
that we have to replace the sum $\sum_{-p 2^{m/2} +1\leq j\leq p
2^{m/2}}$ in $A_{m,n,p}^{(2r)}(t)$, $B_{m,n,p}^{(2r)}(t)$ and
$C_{m,n,p}^{(2r)}(t)$ by $\sum_{-p 2^{\vfrac{m+1}{2}} +1\leq j\leq p
2^{\vfrac{m+1}{2}}}$. And instead of $\sum_{i\geq p 2^{n/2}}$ and $\sum_{i < -p2^{n/2}}$ in $D_{m,n,p}^{(2r)}(t)$,
we must consider $\sum_{i\geq p 2^{\vfrac{n+1}{2}}}$ and $\sum_{i < -p2^{\vfrac{n+1}{2}}}$
respectively. The analysis can then be done mutatis mutandis.

Let us go back to our proof. First, we will prove that
$A_{m,n,p}^{(2r)}(t)$, $B_{m,n,p}^{(2r)}(t)$ and $D_{m,n,p}^{(2r)}(t)$
converge to $0$ in $L^2$ by letting $n$, then $m$, then $p$ tend to
infinity. Second, we will study the convergence in law (in the sense
f.d.d.) of
%
%e3.9 #&#
\begin{equation}
\label{equivalent conv} \bigl( C_{m,n,p}^{(2r)}, 1\leq r \leq N \bigr),
\end{equation}
which will then be equivalent to the convergence in law (in the sense
f.d.d.) of
\[
\bigl( V_n^{(2r)}, 1 \leq r \leq N \bigr).
\]
We will prove that $E[( A_{m,n,p}^{(2r)}(t))^2] \to0$ as $n \to\infty
$. We have, with $\rho$ given by (\ref{rho}) (note that ${\sum_{a \in\Z
} |\rho(a)|^{2r} < \infty}$ if and only if $H < 1 -1/(4r)$, which is
satisfied for any $r \geq1$ because $H < 3/4$),
\begin{eqnarray*}
&&\hspace*{-2pt}E\bigl[\bigl( A_{m,n,p}^{(2r)}(t)\bigr)^2\bigr]
\\
&&\hspace*{-2pt}\quad = \Biggl|2^{-\sfrac{n}{2}}\\
&&\hspace*{-2pt}\qquad{} \times\sum_{-p 2^{m/2}+1\leq j\leq p 2^{m/2}}\sum
_{-p 2^{m/2}+1\leq j'\leq p 2^{m/2}}\sum_{i=(j-1)2^{\vfrac
{n-m}{2}}}^{j2^{\vfrac{n-m}{2}}-1}
\sum_{i'=(j'-1)2^{\vfrac
{n-m}{2}}}^{j'2^{\vfrac{n-m}{2}}-1} E \bigl[H_{2r}
\bigl(X_{i+1}^{(n)}-X_i^{(n)}\bigr)
\\
&&\hspace*{-2pt}\qquad {} \times H_{2r}\bigl(X_{i'+1}^{(n)}-X_{i'}^{(n)}
\bigr) \bigr]E \bigl[\bigl(\mathcal {L}_{i,n}(t)-L_t^{i2^{-n/2}}(Y)
\bigr) \bigl(\mathcal {L}_{i',n}(t)-L_{t}^{i'2^{-n/2}}(Y)
\bigr) \bigr] \Biggr|
\\
&&\hspace*{-2pt}\quad  \leq (2r)! 2^{-\sfrac{n}{2}}\sum_{-p 2^{m/2}+1\leq j\leq p
2^{m/2}}\sum
_{-p 2^{m/2}+1\leq j'\leq p 2^{m/2}}\sum_{i=(j-1)2^{\vfrac
{n-m}{2}}}^{j2^{\vfrac{n-m}{2}}-1}
\sum_{i'=(j'-1)2^{\vfrac
{n-m}{2}}}^{j'2^{\vfrac{n-m}{2}}-1} \rho\bigl(i-i'
\bigr)^{2r}
\\
&&\hspace*{-2pt}\qquad  {}\times\bigl\|\mathcal{L}_{i,n}(t)-L_t^{i2^{-n/2}}(Y)
\bigr\|_2 \times\bigl\|\mathcal {L}_{i',n}(t)-L_t^{i'2^{-n/2}}(Y)
\bigr\|_2,
\end{eqnarray*}
where, in the first equality, we used the independence between $X$ and
$Y$. By the point $3$ of Proposition~\ref{properties-kl}, we have
%
%e3.10 #&#
\begin{equation}
\label{eq:6} \bigl\|\mathcal{L}_{i,n}(t)-L_{t}^{i2^{-n/2}}(Y)
\bigr\|_2 \leq2n2^{-n/4}\|K\| _4\bigl\|L_{t}^{i2^{-n/2}}(Y)
\bigr\|_2^{1/2}.
\end{equation}
On the other hand
%
%e3.11 #&#
\begin{equation}
\label{eq:7} \bigl\|L_{t}^{i2^{-n/2}}(Y)\bigr\|_2 \leq
\bigl\|L_{t}^{i2^{-n/2}}(Y) -L_{t}^0(Y)
\bigr\|_2 + \bigl\|L_{t}^0(Y)\bigr\|_2.
\end{equation}
By the point $2$ of Proposition~\ref{properties-kl}, we have
%
%e3.12 #&#
\begin{equation}
\label{eq:8} \bigl\|L_{t}^{i2^{-n/2}}(Y) -L_{t}^0(Y)
\bigr\|_2 \leq\mu\sqrt{|i|2^{-n/2}} t^{\sfrac{1}{4}}.
\end{equation}
By combining (\ref{eq:7}) and (\ref{eq:8}), we get that ${\|
L_{t}^{i2^{-n/2}}(Y)\|_2 \leq\mu\sqrt{|i|} 2^{-n/4}t^{\sfrac{1}{4}} +
\|L_{t}^0(Y)\|_2} $.

 Since $\sqrt{a+b} \leq\sqrt{a} +\sqrt{b}$ for
all $ a, b \geq0$, we deduce that
%
%e3.13 #&#
\begin{equation}
\label{eq:9} \bigl\|L_{t}^{i2^{-n/2}}(Y)\bigr\|_2^{1/2}
\leq\sqrt{\mu }|i|^{1/4}2^{-n/8}t^{1/8} +
\bigl\|L_{t}^0(Y)\bigr\|_2^{1/2}.
\end{equation}
Finally, $(\ref{eq:9})$ together with $(\ref{eq:6})$ show that
%
%e3.14 #&#
\begin{eqnarray}
\bigl\|\mathcal{L}_{(i,n)}(t)-L_{t}^{i2^{-n/2}}(Y)
\bigr\|_2 &\leq&2\sqrt{\mu}\|K\| _4 t^{1/8}
n2^{-n/4} 2^{-n/8}|i|^{1/4} \nonumber\\[-8pt]\\[-8pt]
&&{}+ 2\|K\|_{4}
\bigl\|L_{t}^0(Y)\bigr\| _2^{1/2}n2^{-n/4}.\nonumber
\end{eqnarray}
As a result,
%
%e3.15 #&#
%e3.16 #&#
%e3.17 #&#
%e3.18 #&#
\begin{eqnarray}
\label{eq:10}&&E\bigl[\bigl( A_{m,n,p}^{(2r)}(t)\bigr)^2\bigr]
\nonumber
\\
&&\quad \leq 4(2r)! \mu t^{1/8}t^{1/8}\|K\|_4^2
2^{-n}2^{-n/4}n^2 \\
&&\qquad {}\times\sum
_{-p
2^{m/2}+1\leq j\leq p 2^{m/2}}\sum_{-p 2^{m/2}+1\leq j'\leq p
2^{m/2}}
 \sum_{i=(j-1)2^{\vfrac{n-m}{2}}}^{j2^{\vfrac{n-m}{2}}-1}\sum
_{i'=(j'-1)2^{\vfrac{n-m}{2}}}^{j'2^{\vfrac{n-m}{2}}-1} \rho \bigl(i-i'
\bigr)^{2r}\bigl|ii'\bigr|^{1/4}
\nonumber
\\
\label
{eq:11}&&\qquad {}+ 4(2r)! \sqrt{\mu}t^{1/8}\|K\|_4^2
\bigl\|L_{t}^0(Y)\bigr\| _2^{1/2}2^{-n}2^{-n/8}n^2\nonumber\\[-8pt]\\[-8pt]
&&\qquad {}\times\sum_{-p 2^{m/2}+1\leq j\leq p 2^{m/2}}
 \sum_{-p 2^{m/2}+1\leq j'\leq p 2^{m/2}}\sum_{i=(j-1)2^{\vfrac
{n-m}{2}}}^{j2^{\vfrac{n-m}{2}}-1}
\sum_{i'=(j'-1)2^{\vfrac
{n-m}{2}}}^{j'2^{\vfrac{n-m}{2}}-1}\rho\bigl(i-i'
\bigr)^{2r}|i|^{1/4}
\nonumber
\\
\label
{eq:12}&&\qquad {}+ 4(2r)! \sqrt{\mu}t^{1/8}\|K\|_4^2
\bigl\|L_t^0(Y)\bigr\| _2^{1/2}2^{-n}2^{-n/8}n^2\nonumber\\[-8pt]\\[-8pt]
&&\qquad {}\times\sum_{-p 2^{m/2}+1\leq j\leq p 2^{m/2}}
 \sum_{-p 2^{m/2}+1\leq j'\leq p 2^{m/2}}\sum_{i=(j-1)2^{\vfrac
{n-m}{2}}}^{j2^{\vfrac{n-m}{2}}-1}
\sum_{i'=(j'-1)2^{\vfrac
{n-m}{2}}}^{j'2^{\vfrac{n-m}{2}}-1}\rho\bigl(i-i'
\bigr)^{2r}\bigl|i'\bigr|^{1/4}
\nonumber
\\
\label{eq:13}&&\qquad {}+ 4(2r)! \|K\|_4^2\bigl\|L_{t}^0(Y)
\bigr\|_2^{1/2}\bigl\|L_{t}^0(Y)\bigr\|
_2^{1/2}2^{-n}n^2\nonumber\\[-8pt]\\[-8pt]
&&\qquad {}\times\sum
_{-p 2^{m/2}+1\leq j\leq p 2^{m/2}}
 \sum_{-p 2^{m/2}+1\leq j'\leq p 2^{m/2}}\sum_{i=(j-1)2^{\vfrac
{n-m}{2}}}^{j2^{\vfrac{n-m}{2}}-1}
\sum_{i'=(j'-1)2^{\vfrac
{n-m}{2}}}^{j'2^{\vfrac{n-m}{2}}-1}\rho\bigl(i-i'
\bigr)^{2r},
\nonumber
\end{eqnarray}
and we are thus left to prove the convergence to $0$ of (\ref
{eq:10})--(\ref{eq:13}) as $n \rightarrow\infty$. Let us do it.
\begin{longlist}[(b)]
\item[(a)] We have
\begin{eqnarray*}
&&2^{-n}n^2\sum_{-p 2^{m/2}+1\leq j\leq p 2^{m/2}}\sum
_{-p
2^{m/2}+1\leq j'\leq p 2^{m/2}}\sum_{i=(j-1)2^{\vfrac
{n-m}{2}}}^{j2^{\vfrac{n-m}{2}}-1}
\sum_{i'=(j'-1)2^{\vfrac
{n-m}{2}}}^{j'2^{\vfrac{n-m}{2}}-1}\rho\bigl(i-i'
\bigr)^{2r}
\\
&&\quad =2^{-n}n^2\sum_{i= -p2^{n/2}}^{p2^{n/2}-1}
\sum_{i'= -p2^{n/2}}^{p2^{n/2}-1} \rho\bigl(i-i'
\bigr)^{2r}
\\
&&\quad \leq2^{-n}n^2\sum_{i= -p2^{n/2}}^{p2^{n/2}-1}
\sum_{i'\in\Z}\rho \bigl(i'
\bigr)^{2r}=\sum_{i'\in\Z}\rho
\bigl(i'\bigr)^{2r} n^22^{-n}
\bigl(2p2^{n/2}\bigr).
\end{eqnarray*}
Since it is clear that the last quantity converges to $0$ as $n \to
\infty$, one deduces that $(\ref{eq:13})$ tends to zero.
\item[(b)] Since $ -p 2^{m/2}+1\leq j' \leq p 2^{m/2}$ and
$(j'-1)2^{\vfrac{n-m}{2}}\leq i'\leq j'2^{\vfrac{n-m}{2}}-1$, we deduce
that $-p 2^{n/2} \leq i'\leq p 2^{n/2} -1$. So, $|i'| \leq p2^{n/2}$.
Consequently we have that $|i'|^{1/4} \leq p^{1/4}2^{n/8}$, which shows that
\begin{eqnarray*}
&&\hspace*{-4pt}2^{-n}2^{-n/8}n^2\\
&&\hspace*{-7pt}\qquad \times\sum
_{-p 2^{m/2}+1\leq j\leq p 2^{m/2}}\sum_{-p
2^{m/2}+1\leq j'\leq p 2^{m/2}}\sum
_{i=(j-1)2^{\vfrac
{n-m}{2}}}^{j2^{\vfrac{n-m}{2}}-1} \sum_{i'=(j'-1)2^{\vfrac
{n-m}{2}}}^{j'2^{\vfrac{n-m}{2}}-1}
\rho\bigl(i-i'\bigr)^{2r}\bigl|i'\bigr|^{1/4}
\\
&&\hspace*{-6pt}\quad \leq p^{1/4}2^{-n}n^2\\
&&\hspace*{-7pt}\qquad \times\sum
_{-p 2^{m/2}+1\leq j\leq p 2^{m/2}}\sum_{-p
2^{m/2}+1\leq j'\leq p 2^{m/2}}\sum
_{i=(j-1)2^{\vfrac
{n-m}{2}}}^{j2^{\vfrac{n-m}{2}}-1} \sum_{i'=(j'-1)2^{\vfrac
{n-m}{2}}}^{j'2^{\vfrac{n-m}{2}}-1}
 \rho\bigl(i-i'\bigr)^{2r}
\end{eqnarray*}
and this last quantity converges to $0$ by the same argument as above.
This shows that (\ref{eq:12}) tends to zero.
\item[(c)] Following the same strategy as in point (b), one deduces
that (\ref{eq:11}) tends to zero. Details are left to the reader.
\item[(d)] By the same arguments as above, one can see that
$|ii'|^{1/4} \leq p^{1/2}2^{n/4}$. It follows that
\begin{eqnarray*}
&&\hspace*{-4pt}2^{-n}2^{-n/4}n^2 \\
&&\hspace*{-11pt}\qquad \times\sum
_{-p 2^{m/2}+1\leq j\leq p 2^{m/2}}\sum_{-p
2^{m/2}+1\leq j'\leq p 2^{m/2}}\sum
_{i=(j-1)2^{\vfrac
{n-m}{2}}}^{j2^{\vfrac{n-m}{2}}-1} \sum_{i'=(j'-1)2^{\vfrac
{n-m}{2}}}^{j'2^{\vfrac{n-m}{2}}-1}
 \rho\bigl(i-i'\bigr)^{2r}\bigl|ii'\bigr|^{1/4}
\\
&&\hspace*{-7pt}\quad  \leq p^{1/2}2^{-n}n^2\\
&&\hspace*{-11pt}\qquad \times \sum
_{-p 2^{m/2}+1\leq j\leq p 2^{m/2}}\sum_{-p
2^{m/2}+1\leq j'\leq p 2^{m/2}}\sum
_{i=(j-1)2^{\vfrac
{n-m}{2}}}^{j2^{\vfrac{n-m}{2}}-1} \sum_{i'=(j'-1)2^{\vfrac
{n-m}{2}}}^{j'2^{\vfrac{n-m}{2}}-1}
 \rho\bigl(i-i'\bigr)^{2r},
\end{eqnarray*}
which converges to $0$ by the same arguments as above. Hence, (\ref
{eq:10}) tends to zero. The proof of $E[( A_{m,n,p}^{(2r)}(t))^2] \to
0$ as $ n\to\infty$ is complete.
\end{longlist}

Now, let us prove the convergence of $B_{m,n,p}^{(2r)}(t)$ to $0$ in
$L^2$ as $m \to \infty$, uniformly in $n$. We have
\begin{eqnarray*}
&&E\bigl[\bigl( B_{m,n,p}^{(2r)}(t)\bigr)^2\bigr]
\\
&&\quad =2^{-n/2} \sum_{-p 2^{m/2}+1\leq j\leq p 2^{m/2}} \sum
_{-p 2^{m/2}+1\leq j'\leq
p 2^{m/2}} \sum_{i=(j-1)2^{\vfrac{n-m}{2}}}^{j2^{\vfrac{n-m}{2}}-1}
\sum_{i'=(j'-1)2^{\vfrac{n-m}{2}}}^{j'2^{\vfrac{n-m}{2}}-1}(2r)! \rho\bigl(i-i'\bigr)^{2r}
\\
&&\qquad{}  \times E\bigl[
\bigl(L_t^{i2^{-n/2}}(Y)-L_t^{j2^{-m/2}}(Y)
\bigr) \bigl(L_t^{i'2^{-n/2}}(Y)-L_t^{j'2^{-m/2}}(Y)
\bigr)\bigr].
\end{eqnarray*}
By Proposition~\ref{properties-kl} (point 2) and Cauchy--Schwarz, there
is a universal constant $\mu$ such that
\begin{eqnarray*}
&& \bigl|E\bigl[\bigl(L_t^{i2^{-n/2}}(Y)-L_t^{j2^{-m/2}}(Y)
\bigr) \bigl(L_t^{i'2^{-n/2}}(Y)-L_t^{j'2^{-m/2}}(Y)
\bigr)\bigr] \bigr|
\\
&&\quad  \leq \mu^2 \sqrt{t} \sqrt {\bigl|i2^{-n/2}-j2^{-m/2}\bigr|\bigl|i'2^{-n/2}-j'2^{-m/2}\bigr|}
\leq\mu^2 \sqrt{t}2^{-m/2}.
\end{eqnarray*}
This yields
\begin{eqnarray*}
\sup_n E\bigl[\bigl( B_{m,n,p}^{(2r)}(t)
\bigr)^2\bigr] &\leq&\mu^2 (2r)! 2^{-m/2}\sqrt{t}
\\
&&{} \times\sup_n \Biggl\{2^{-n/2} \sum
_{i= -p2^{n/2}}^{p2^{n/2}-1} \sum_{i'= -p2^{n/2}}^{p2^{n/2}-1}
\rho\bigl(i-i'\bigr)^{2r} \Biggr\}
\\
&\leq&\mu^2 (2r)! 2^{-m/2}\sqrt{t} 2p\sum
_{i\in\mathbb{Z}} \rho(i)^{2r} ,
\end{eqnarray*}
which converges to $0$ as $m \to\infty$.

Finally, let us prove that $D_{m,n,p}^{(2r)}(t)$ converges to $0$ in
$L^2$ as $p \to\infty$, uniformly in $m$ and $n$. We have
%
%e3.19 #&#
%e3.20 #&#
%e3.21 #&#
\begin{eqnarray}
\label{eq:correction1}&&\hspace*{-35pt}E \bigl[ \bigl(D_{m,n,p}^{(2r)}(t)\bigr)^2 \bigr]
\nonumber
\\
&&\hspace*{-35pt}\quad = 2^{-\sfrac{n}{2}}\sum_{i \geq p2^{n/2}}\sum
_{j \geq p2^{n/2}}E \bigl[H_{2r}\bigl(X_{i+1}^{(n)}
- X_i^{(n)}\bigr)H_{2r}\bigl(X_{j+1}^{(n)}
- X_j^{(n)}\bigr) \mathcal{L}_{i,n}(t)
\mathcal{L}_{j,n}(t) \bigr]
\\
\label{eq:correction2}&&\hspace*{-35pt} \qquad {}+ 22^{-\sfrac{n}{2}}\sum_{i \geq p2^{n/2}}\sum
_{j < -p2^{n/2}}E \bigl[H_{2r}\bigl(X_{i+1}^{(n)}
- X_i^{(n)}\bigr)H_{2r}\bigl(X_{j+1}^{(n)}
- X_j^{(n)}\bigr) \mathcal{L}_{i,n}(t)\mathcal{L}_{j,n}(t) \bigr]
\\
\label{eq:correction3}&&\hspace*{-35pt}\qquad {} + 2^{-\sfrac{n}{2}}\sum_{i < -p2^{n/2}}\sum
_{j < -p2^{n/2}}E \bigl[H_{2r}\bigl(X_{i+1}^{(n)}
- X_i^{(n)}\bigr)H_{2r}\bigl(X_{j+1}^{(n)}
- X_j^{(n)}\bigr)
\mathcal{L}_{i,n}(t)\mathcal{L}_{j,n}(t) \bigr],\nonumber\\[-6pt]
\end{eqnarray}
and we are thus left to prove the convergence to $0$ of (\ref
{eq:correction1})--(\ref{eq:correction3}) as $p \to\infty$, uniformly
in $m$ and $n$. Let us do it.
\begin{longlist}[(b)]
\item[(a)] We have
\begin{eqnarray*}
&& \biggl|2^{-\sfrac{n}{2}}\sum_{i \geq p2^{n/2}}\sum
_{j \geq p2^{n/2}}E \bigl[H_{2r}\bigl(X_{i+1}^{(n)}
- X_i^{(n)}\bigr)H_{2r}\bigl(X_{j+1}^{(n)}
- X_j^{(n)}\bigr) \mathcal{L}_{i,n}(t)
\mathcal{L}_{j,n}(t) \bigr] \biggr|
\\
&&\quad = \biggl|2^{-\sfrac{n}{2}}\sum_{i \geq p2^{n/2}}\sum
_{j \geq p2^{n/2}}E \bigl[H_{2r}\bigl(X_{i+1}^{(n)}
- X_i^{(n)}\bigr)H_{2r}\bigl(X_{j+1}^{(n)}
- X_j^{(n)}\bigr) \bigr]
 E \bigl[\mathcal{L}_{i,n}(t)\mathcal{L}_{j,n}(t)
\bigr] \biggr|
\\
&&\quad = (2r)! 2^{-\sfrac{n}{2}} \sum_{i \geq p2^{n/2}}\sum
_{j \geq
p2^{n/2}}\rho(i-j)^{2r}E \bigl[ \mathcal{L}_{i,n}(t)
\mathcal {L}_{j,n}(t) \bigr],
\end{eqnarray*}
where, in the second equality, we used the independence between $X$ and
$Y$. It is enough to prove that, uniformly in $n$ and $m$, and as $p
\to\infty$:
%
%e3.22 #&#
\begin{equation}
\label{eq:14} 2^{-\sfrac{n}{2}} \sum_{i \geq p2^{n/2}}\sum
_{j \geq p2^{n/2}}\rho (i-j)^{2r}E \bigl[
\mathcal{L}_{i,n}(t)\mathcal{L}_{j,n}(t) \bigr] \to0.
\end{equation}
We can write
\begin{eqnarray*}\label{eq:15}
&&2^{-\sfrac{n}{2}} \sum_{i \geq p2^{n/2}}\sum
_{j \geq p2^{n/2}}\rho (i-j)^{2r}E \bigl[ \mathcal{L}_{i,n}(t)
\mathcal{L}_{j,n}(t) \bigr]
\\
&&\quad \leq2^{-\sfrac{n}{2}}\sum_{i\geq p2^{n/2}}\sum
_{j \geq p2^{n/2}}\rho (i-j)^{2r}E \biggl[ \frac{1}{2}
\bigl( \mathcal{L}_{i,n}(t)^2 + \mathcal
{L}_{j,n}(t)^2 \bigr) \biggr]
\\
&&\quad = 2^{-\sfrac{n}{2}}\sum_{i \geq p2^{n/2}}E \bigl[
\mathcal{L}_{i,n}(t)^2 \bigr]\sum
_{j \geq p2^{n/2}}\rho(i-j)^{2r}
\\
&&\quad  \leq2^{-\sfrac{n}{2}}\sum_{i \geq p2^{n/2}}E \bigl[ \mathcal
{L}_{i,n}(t)^2 \bigr]\sum_{j \in\Z}
\rho(j)^{2r} =C_r 2^{-\sfrac
{n}{2}}\sum
_{i \geq p2^{n/2}}E \bigl[ \mathcal{L}_{i,n}(t)^2
\bigr],
\end{eqnarray*}
where ${C_r:=\sum_{j \in\Z}\rho(j)^{2r} < \infty}$.
By the third point of Proposition~\ref{properties-kl}, we have
\[
\bigl|\mathcal{L}_{i,n}(t)\bigr| \leq L_t^{i2^{-n/2}}(Y)+2Kn2^{-n/4}
\sqrt {L_t^{i2^{-n/2}}(Y)}
\]
so that
%
%e3.23 #&#
\begin{equation}
\label{eq:16} E \bigl[ \mathcal{L}_{i,n}(t)^2 \bigr]
\leq2E \bigl[ L_t^{i2^{-n/2}}(Y)^2 \bigr] +
8n^2 2^{-n/2}\bigl\|K^2\bigr\|_2
\bigl\|L_t^{i2^{-n/2}}(Y)\bigr\|_2.
\end{equation}
On the other hand, thanks to the point 1 of Proposition~\ref
{properties-kl}, we have
%
%e3.24 #&#
\begin{equation}
\label{eq:17} E \bigl[ L_t^{i2^{-n/2}}(Y)^2 \bigr]
\leq C t \exp \biggl(-\frac
{(i2^{-n/2})^2}{2t} \biggr).
\end{equation}
Consequently, we get
%
%e3.25 #&#
\begin{equation}
\label{eq:18} \bigl\|L_{t_{a_2}}^{i2^{-n/2}}(Y)\bigr\|_2 \leq
C^{1/2} t^{1/2} \exp \biggl(-\frac
{(i2^{-n/2})^2}{4t} \biggr).
\end{equation}
By combining (\ref{eq:16}) with (\ref{eq:17}) and (\ref{eq:18}), we
deduce that
\begin{eqnarray*}
2^{-\sfrac{n}{2}}\sum_{i \geq p2^{n/2}}E \bigl[
\mathcal{L}_{i,n}(t)^2 \bigr]&\leq&2C t 2^{-n/2}\sum
_{i \geq p2^{n/2}}\exp \biggl(-\frac
{(i2^{-n/2})^2}{2t} \biggr)
\\
&&{} + 8C^{1/2}t^{1/2}\bigl\|K^2\bigr\|_2
n^2 2^{-n/2}
\\
&&\hphantom{{}+}{} \times2^{-n/2}\sum_{i \geq p2^{n/2}}\exp \biggl(-
\frac
{(i2^{-n/2})^2}{4t} \biggr).
\end{eqnarray*}
But, for $a \in\{2,4\}$,
\[
{2^{-n/2}\sum_{i \geq p2^{n/2}}\exp \biggl(-
\frac{(i2^{-n/2})^2}{at} \biggr) \leq \int_{p-1}^{\infty}\exp
\biggl(\frac{-x^2}{at} \biggr)\, \mathrm{d}x \displaystyle \mathop{\longrightarrow}_{p \rightarrow\infty}
0}.
\]
This proves (\ref{eq:14}). Hence, we deduce that (\ref{eq:correction1})
converges to $0$ as $p \rightarrow\infty$ uniformly in $n$ and~$m$.
\item[(b)]Following the same strategy as in point (a), one deduces
that (\ref{eq:correction2}) and (\ref{eq:correction3}) converge to $0$
as $p \rightarrow\infty$ uniformly in $n$ and $m$. Details are left to
the reader.
\end{longlist}
This shows that $D_{m,n,p}^{(2r)}(t)$ converges to $0$ in $L^2$ as $p
\to\infty$, uniformly in $m$ and $n$.

To finish our proof of (\ref{conv2}), it remains to prove that, by
letting $n$, then $m$, then $p$ tend to infinity, we get
%
%e3.26 #&#
\begin{equation}
\label{proof-conv2} \bigl\{ C_{m,n,p}^{(2r)}(t), 1\leq r \leq N \bigr
\}_{t \geq0} \stackrel{\mathrm{f.d.d.}} {\longrightarrow} \biggl\{
\sigma_{2r} \int_{-\infty}^\infty
L_t^x(Y)\,\mathrm{d}B^{(r)}_x \dvt  1\leq r\leq N
\biggr\}_{t\geq0}.
\end{equation}
Since $H < 3/4$, we claim that, as $n \rightarrow \infty$,
%
%e3.27 #&#
\begin{eqnarray}
\label{eq:18bis}
&&\Biggl(2^{-n/4}\sum_{i=(j-1)2^{\vfrac{n-m}{2}}}^{j2^{\vfrac{n-m}{2}}-1}
H_{2r}\bigl(X_{i+1}^{(n)} - X_i^{(n)}
\bigr), 1 \leq r \leq N \dvt  -p2^{m/2}+1 \leq j \leq p2^{m/2}
\Biggr) \nonumber
\\[-8pt]\\[-8pt]
&&\quad \stackrel{\mathrm{law}} {\longrightarrow}\bigl( \sigma_{2r}\bigl(B_{(j+1)2^{-m/2}}^{(r)} -
B_{j2^{-m/2}}^{(r)}\bigr), 1 \leq r \leq N \dvt  -p2^{m/2}+1
\leq j \leq p2^{m/2} \bigr),\nonumber
\end{eqnarray}
where $(B^{(1)}, \ldots, B^{(N)})$ is a $N$-dimensional two-sided
Brownian motion.

Indeed, it is clear, using the self-similarity property of $X$, that
the convergence in law of
\[
\Biggl(2^{-n/4}\sum_{i=(j-1)2^{\vfrac{n-m}{2}}}^{j2^{\vfrac{n-m}{2}}-1}
H_{2r}\bigl(X_{i+1}^{(n)} - X_i^{(n)}
\bigr), 1 \leq r \leq N \dvt  -p2^{m/2}+1 \leq j \leq p2^{m/2}
\Biggr)
\]
is equivalent to the convergence in law of
\[
\Biggl(2^{-n/4}\sum_{i=(j-1)2^{\vfrac{n-m}{2}}}^{j2^{\vfrac{n-m}{2}}-1}
H_{2r}(X_{i+1} - X_i), 1 \leq r \leq N \dvt
-p2^{m/2}+1 \leq j \leq p2^{m/2} \Biggr).
\]
Then, Breuer--Major Theorem \ref{th:bm} applies and yields that, as $n\to
\infty$ and for any fixed
 $1 \leq r \leq N$,
\begin{eqnarray*}
&&\Biggl(2^{-n/4}\sum_{i=(j-1)2^{\vfrac{n-m}{2}}}^{j2^{\vfrac{n-m}{2}}-1}
H_{2r}(X_{i+1} - X_i) \dvt  -p2^{m/2}+1
\leq j \leq p2^{m/2} \Biggr)
\\
&&\quad  \stackrel{\mathrm{law}} {\longrightarrow}\bigl( \sigma_{2r}\bigl(B_{(j+1)2^{-m/2}}^{(r)} -
B_{j2^{-m/2}}^{(r)}\bigr)\dvt  -p2^{m/2}+1 \leq j \leq
p2^{m/2} \bigr).
\end{eqnarray*}
In addition, from Peccati--Tudor Theorem \ref{th:pec} and taking into
account the orthogonality of Hermite polynomial with different orders,
we deduce (\ref{eq:18bis}). (The detailed proof of this result is
similar to the proof of (\ref{cv2}).)

As a consequence of (\ref{eq:18bis}), and thanks to the independence
of $X$ and $Y$, we have that as $n \to \infty$,
\begin{eqnarray*}
&& \bigl\{ C_{m,n,p}^{(2r)}(t), 1\leq r \leq N \bigr
\}_{t \geq0}
\\
&&\quad  \stackrel{\mathrm{f.d.d.}} {\longrightarrow}\Biggl\{ \sigma_{2r}\sum_{j=-p2^{m/2}+1}^{p2^{m/2}}L_t^{j2^{-m/2}}(Y)
\bigl(B_{(j+1)2^{-m/2}}^{(r)} - B_{j2^{-m/2}}^{(r)}\bigr),
1\leq r \leq N \Biggr\}_{t \geq0}.
\end{eqnarray*}
Since, for any fixed $t \geq0$ and $1\leq r \leq N$ and as $m \to
\infty$,
\[
\sum_{j=-p2^{m/2}+1}^{p2^{m/2}}L_t^{j2^{-m/2}}(Y)
\bigl(B_{(j+1)2^{-m/2}}^{(r)} - B_{j2^{-m/2}}^{(r)}\bigr)
\stackrel{P} {\longrightarrow}\int_{-p}^{p}L_t^x(Y)
\,\mathrm{d}B_x^{(r)},
\]
and since ${\int_{-p}^{p}L_t^x(Y) \,\mathrm{d}B_x^{(r)} \stackrel
{P}{\longrightarrow}\int_{\R}L_t^x(Y)\,\mathrm{d}B_x^{(r)}}$ as $p \to \infty$,
we deduce finally that by letting $m$, then $p$ tend to infinity, we get
\begin{eqnarray*}
&& \Biggl\{ \sigma_{2r}\sum_{j=-p2^{m/2}+1}^{p2^{m/2}}L_t^{j2^{-m/2}}(Y)
\bigl(B_{(j+1)2^{-m/2}}^{(r)} - B_{j2^{-m/2}}^{(r)}\bigr),
1\leq r \leq N \Biggr\}_{t \geq0}
\\
&&\quad  \stackrel {\mathrm{f.d.d.}} {
\longrightarrow}\biggl\{ \sigma_{2r} \int_{-\infty}^\infty
L_t^x(Y)\,\mathrm{d}B^{(r)}_x \dvt  1\leq r\leq N
\biggr\}_{t\geq0}.
\end{eqnarray*}
This proves (\ref{proof-conv2}), and consequently (\ref{conv2}).

%s3.3 #&#
\subsection{Proof of Corollary \texorpdfstring{\protect\ref{cor-main}}{1.2}}\label{sec3.3}
Let us decompose $x^p$ in terms of Hermite polynomials. We have
${x^p= \sum_{k=0}^p a_{p,k}H_k(x)}$, where $a_{p,k}$ is some (explicit)
integer. To calculate $a_{p,k}$, let $N$ be a centred Gaussian variable
with variance one. We have
%
%e3.28 #&#
\begin{equation}
\label{eq:19} N^p= \sum_{k=0}^p
a_{p,k}H_k(N).
\end{equation}
Thanks to the orthogonality property of Hermite polynomials with
different orders and to the well known fact that $E[H_k(N)^2]=k!$, we get
%
%e3.29 #&#
\begin{equation}
\label{eq:20} a_{p,k} = \frac{1}{k!} E\bigl[N^pH_k(N)
\bigr].
\end{equation}
On the other hand (see, e.g., \cite{8}, page 19) we have, for all $k
\geq1$,
%
%e3.30 #&#
\begin{equation}
\label{eq:21} H_k(x)= \sum_{l=0}^{\lfloor k/2 \rfloor}
\frac{k!(-1)^l}{l!(k-2l)!2^l}x^{k-2l}.
\end{equation}
By combining (\ref{eq:20}) with (\ref{eq:21}), we deduce that
%
%e3.31 #&#
\begin{equation}
\label{eq:22} a_{p,k}=\sum_{l=0}^{\lfloor k/2 \rfloor}
\frac
{(-1)^l}{l!(k-2l)!2^l}E\bigl(N^{p+k-2l}\bigr).
\end{equation}
Thus,
\[
a_{p,k}=\lleft\{ %
\begin{array} {l@{\qquad}l} \displaystyle \sum
_{l=0}^{\lfloor k/2 \rfloor}\displaystyle \frac
{(-1)^l(p+k-2l)!}{l!(k-2l)!2^l2^{\vfrac{p+k-2l}{2}}(\vfrac{p+k-2l}{2})!} & \mbox{if $p$ and
$k$ are odd},
\\\noalign{\vspace*{2pt}}
\displaystyle \sum_{l=0}^{\lfloor k/2 \rfloor}\displaystyle \frac
{(-1)^l(p+k-2l)!}{l!(k-2l)!2^l2^{\vfrac{p+k-2l}{2}}(\vfrac{p+k-2l}{2})!} & \mbox{if $p$ and $k$ are even},
\\\noalign{\vspace*{2pt}}
0 & \mbox{otherwise}. \end{array} %
 \rright.
\]
As a result, we deduce that if $p$ is odd, then
%
%e3.32 #&#
\begin{equation}
\label{eq:23} x^p=\sum_{k=1}^{\lfloor p/2 \rfloor+1}
a_{p,2k-1}H_{2k-1}(x),
\end{equation}
whereas if $p$ is even, then
%
%e3.33 #&#
\begin{equation}
\label{eq:24} x^p=\sum_{k=0}^{p/2}
a_{p,2k}H_{2k}(x),
\end{equation}
Finally, thanks to (\ref{eq:23}), (\ref{eq:24}), Theorem~\ref{main} and
the Continuous Mapping theorem, we deduce the content of Corollary~\ref
{cor-main}.

% zodis "Acknowledgments" paliekamas pagal autoriu
\section*{Acknowledgements}
I thank my supervisor Ivan Nourdin for his help and for several
interesting discussions about this article.
Also, I thank an anonymous referee for his/her many valuable comments
and remarks on a previous version of this work.

%suskaldyti doi

% imsref loaded by jurgita.kaciuliene, 2014-02-28 10:09:26
% imsref loaded by jurgita.kaciuliene, 2014-02-28 10:10:38

\printhistory

\end{document}